# How to use the Kohonen algorithm to simultaneously analyze individuals and modalities in a survey

Marie Cottrell, Patrick Letrémy,

Université Paris I, SAMOS-MATISSE, UMR CNRS 8595
90 rue de Tolbiac,
F-75634 Paris Cedex 13, France
pley,cottrell@univ-paris1.fr

**Abstract :**
The Kohonen algorithm (SOM, Kohonen,1984, 1995) is a very powerful tool for data analysis. It was originally designed to model organized connections between some biological neural networks. It was also immediately considered as a very good algorithm to realize vectorial quantization, and at the same time pertinent classification, with nice properties for visualization. If the individuals are described by quantitative variables (ratios, frequencies, measurements, amounts, etc.), the straightforward application of the original algorithm leads to build code vectors and to associate to each of them the class of all the individuals which are more similar to this code-vector than to the others. But, in case of individuals described by categorical (qualitative) variables having a finite number of modalities (like in a survey), it is necessary to define a specific algorithm. In this paper, we present a new algorithm inspired by the SOM algorithm, which provides a simultaneous classification of the individuals and of their modalities.

## 1. Introduction

For real-valued databases, many tools are available (descriptive techniques, factorial data analysis, neural networks, Kohonen maps, prediction models, and so on).

Recently a lot of methods have been proposed to handle discrete databases, mainly in text analysis or gene expression analysis frameworks. The authors have defined discrete analogues to PCA, including probabilistic latent semantic analysis (Hofmann, 1999), latent Dirichlet allocation (Blei, Ng & Jordan, 2002), generative aspect models (Minka & Lafferty, 2002), Multinomial PCA (Buntime & Perttu, 2003), biclustering algorithms (Cheng & Church). All these methods are adapted to process very high dimensional discrete data. They rely on probabilistic hypotheses, where the data are viewed as a mixture of Gaussians or of multinomial distributions.

In this paper we are interested in databases which come from surveys (for example socio-economics surveys or medical surveys). The individuals are mainly described by categorical (i.e. qualitative) variables. In that case, each individual has to answer a number of questions, each of them having a finite number of possible modalities (i.e. sex, professional category, level of income, kind of employment, place of housing, type of car, level of education, etc.). This way to gather the answers corresponds to a legal obligation, since nobody has to be recognized from his answers. The databases are structured, since there is only one possible answer for each question. After cleaning the data, the total number of modalities is generally medium (about 100 as a maximum).

Most of the time, the modalities are encoded by integer values (1, 2, 3,...), and sometimes are viewed as numerical values. But in fact, it is well known that this is not adequate. The encoding values can be not comparable, the codes are neither necessarily put in order nor regularly spaced (for example, is it blue color less than brown color?, how to arrange in order the types of car, the places of housing,...?). In that case, using the encoding of the modalities as quantitative (real-valued) variables has no meaning and this kind of qualitative databases need a specific treatment.

In this paper, we present a technique that is analogue to the Correspondence Analysis and does not assume any probabilistic hypotheses. It is mainly descriptive and does not rely on likelihood writing. Its visualization properties are fully useful when the number of modalities is not too large

Let us define the notations. We consider a sample of $N$ individuals and a number $K$ of categorical variables (i.e. the questions of the survey). Each variable $k = 1,2, ..., K$ has $m_k$ possible modalities. For each individual, there is one and only one possible modality for each question. So, if $M$ is the total number of modalities, each individual is represented by a row $M$-vector with values in $\{0, 1\}$. There is only one 1 between the 1st component and the $m_1$-th one, only one 1 between the $(m_1+1)$-th component and the $(m_1+m_2)$-th one and so on. The $(N\times M)$ data matrix is called the *complete disjunctive table* and is denoted by $D = (d_{ij})$, $i = 1, .., N$, $j = 1,..., M$. The term $d_{ij}$ takes its values in $\{0,1\}$. This table $D$ contains all the information about the individuals. See in Figure 1, a stylized representation of such a table.

| | $m_1$ | | | $m_2$ | | $m_3$ | | |
|---|---|---|---|---|---|---|---|---|
| *Ind* | 1 | 2 | 3 | 1 | 2 | 1 | 2 | 3 |
| 1 | 0 | 1 | 0 | 0 | 1 | 1 | 0 | 0 |
| 2 | 1 | 0 | 0 | 1 | 0 | 0 | 1 | 0 |
| … | | | | | | | | |
| … | | | | | | | | |
| *i* | 0 | 0 | 1 | 0 | 1 | 0 | 0 | 1 |
| | | | | | | | | |
| | | | | | | | | |
| *N* | 0 | 0 | 1 | 1 | 0 | 0 | 0 | 1 |

*Fig. 1: A Complete Disjunctive Table*

Two kinds of problems can be handled : it can be of interest to only study the relations between modalities. This kind of study put in evidence what are the modalities which are strongly associated and is able to define clusters among them. In that case, it is not necessary to work on the Complete Disjunctive Table. One can summarize the data into a Burt Table which is a Contingency Hypercube, where all variables are crossed two by two. It is then possible to deal with this Burt Table with factorial methods (Burt, 1950, Benzécri, 1992) or Kohonen-inspired methods (Cottrell and Ibbou, 1995, Cottrell et al. 1999).

But it can be more interesting and more valuable to classify at the same time the individuals and their modalities. This is the goal of this paper, that is to represent the individuals near their own modalities, and the modalities near the individuals who share them. The position of the modalities among the classes can be controlled by the distribution of the individuals who belong to the same classes.

To achieve this goal, in classical linear data analysis, one uses a factorial method, the Multiple Correspondence Analysis (MCA) which is a variant of Principal Component Analysis. See Benzécri, 1992, Lebart et al., 1984, for theory and examples.

Let us recall briefly how to achieve a Multiple Correspondence Analysis on the Complete Disjunctive Table. It is very similar to the analysis of a contingency table, where the row-variable would be the INDIVIDUAL variable, and the column-variable would be the MODALITY variable.

One defines successively

- the row-margin of the table: $d_{i.} = \sum_{j=1}^{M} d_{ij} = K$, for every *i*, since it is the number of questions,
- the colum-margin of the table: $d_{.j} = \sum_{i=1}^{N} d_{ij}$, it is the number of individuals who chose modality *j*.
- the row-profiles (which sum to 1), with entry $\frac{d_{ij}}{d_{i.}}$
- the column-profiles (which sum to 1), with entry $\frac{d_{ij}}{d_{.j}}$

So two set of points are defined: *N* row-points in $R^M$ corresponding to the individuals and *M* column-points in $R^N$ corresponding to the modalities.

After that, we can note that the usual Euclidean distance between the rows-profiles gives an important weight to the modalities which are more frequent than the others. To take into account this fact and to correct it, we weight each difference by the inverse of the frequency by defining the distance between rows *i* and *i'* :

$$\chi^2(i,i') = \sum_{j=1}^{M} \frac{1}{d_{.j}} \left( \frac{d_{ij}}{d_{i.}} - \frac{d_{i'j}}{d_{i'.}} \right)^2 .$$

This distance is called Chi-2 distance, since it is one of the term that are added up to compute the Chi-2 statistic which is used to test the independence between the colums and the rows.

This distance between rows can be re-written as :

$$\chi^2(i,i') = \sum_{j=1}^{M} \frac{1}{d_{.j}} \left( \frac{d_{ij}}{d_{i.}} - \frac{d_{i'j}}{d_{i'.}} \right)^2 = \sum_{j=1}^{M} \left( \frac{d_{ij}}{d_{i.}\sqrt{d_{.j}}} - \frac{d_{i'j}}{d_{i'.}\sqrt{d_{.j}}} \right)^2 = \sum_{j=1}^{M} \left( \frac{d_{ij}}{K\sqrt{d_{.j}}} - \frac{d_{i'j}}{K\sqrt{d_{.j}}} \right)^2 .$$

The first correction on Table *D*, is to consider as entry $\frac{d_{ij}}{K\sqrt{d_{.j}}}$ instead of $d_{ij}$.

The same distance is defined between the colums by

$$\chi^2(j,j') = \sum_{i=1}^{N} \frac{1}{d_{i.}} \left( \frac{d_{ij}}{d_{.j}} - \frac{d_{ij'}}{d_{.j'}} \right)^2 = \sum_{i=1}^{N} \left( \frac{d_{ij}}{\sqrt{d_{i.}}d_{.j}} - \frac{d_{i'j}}{\sqrt{d_{i'.}}d_{.j}} \right)^2 = \sum_{i=1}^{N} \left( \frac{d_{ij}}{\sqrt{K}d_{.j}} - \frac{d_{i'j}}{\sqrt{K}d_{.j}} \right)^2 .$$

The corresponding correction consists in replacing $d_{ij}$ by $\frac{d_{ij}}{\sqrt{K}d_{.j}}$.

With these distances, two modalities chosen by the same individuals are coincident, the modalities chosen by few individuals are far from the others. In the same way, two individuals who chose the same modalities are close, and are different when they did not answer in the same way. More a modality is rare, more it is important in the calculus of the distances.

We can see that until this point, the corrections are not the same for rows or for columns.

In the row-space (the individuals), the next step is to use a Principal Component Analysis on the *N* rows, endowed with the Chi-2 distance, and weighted by the frequency *K* for each row). Each squared term has to be multiplied by *K*, so it is equivalent to use a Principal Component Analysis on *N* corrected rows with entry

$$d_{ij}^{c} = \frac{d_{ij}}{\sqrt{K}\sqrt{d_{.j}}}.$$

Let us denote by $D^c$ the Corrected Table whose entries are $d_{ij}^{c}$.

The next step is the computation of the eigenvalues and eigenvectors of the $M \times M$-matrix $D^{c}{}'D^{c}$, put in order according to decreasing values of the eigenvalues. The eigenvectors define the principal axes, and the successive projections on axes 1-2, 1-3, 2-3, and so on, provide successive representations of the rows (i.e. the individuals) with decreasing significance. The number of axes is the rank of the matrix $D^{c}{}'D^{c}$, that is $M - K$, the difference between the number of modalities and the number of questions, since for each question the sum of columns is equal to 1.

As to the colums (the modalities), we can observe that given that the formula which defines $d_{ij}^{c}$ is symmetric with respect to *i* and *j*, the Principal Component Analysis on the *M* columns consists in computing the eigenvalues and the eigenvectors of the $N \times N$-matrix $D^{c}D^{c}{}'$. The number of non-zero eigenvalues of both matrices is the same ($M - K$) and there exist close relations between principal axes of both analysis. So it is possible to superpose the projections of the individuals and of the modalities on each two-dimensional map. The main property is that each modality is drawn as an approximate gravity center of the individuals that possess it, and that each individual is approximately at the gravity center of his own modalities.

But the approximation can be very poor. It is a projection method, which provides several two-dimensional maps, each of them representing a small percent of the global information. So it is necessary to look at several maps at once, the modalities and/or the individuals are more or less well represented, and it is not always easy to deduce pertinent conclusions about the proximity between modalities, between individuals, and between modalities and individuals. A main drawback of the MCA representations is that the linear projection can distort the distances. Two neighbor points are always projected on neighbor locations, but the converse is not true : two apparently neighbor locations in the projection subspace can correspond to very distant points in an orthogonal subspace.

To overcome these drawbacks, we propose as an alternative to the classical factorial method we just described, to use an adaptation of the Kohonen algorithm. It is well known that this algorithm is able to provide nice representation of numerical data, analogous to a projection on the two first axes of a Principal Component Analysis. (Kaski, 1997). The Kohonen maps are easy to interpret since the provided classification respects the topology of the input.

However, as we observed in the beginning of the introduction, in its classical version, the Kohonen algorithm is not adapted to deal with qualitative data. Hence the main point of our work is to propose modifications of the genuine SOM algorithm, designed to deal with individuals described with qualitative data, to simultaneously represent the individuals and the modalities, by classifying and visualizing them on a Kohonen map.

## 2. A new algorithm for a simultaneous analysis of individuals together the modalities

In the same way as explained before, we consider the Complete Disjunctive Table *D* as a contingency table, which crosses an extra-variable MODALITY with an extra-variable INDIVIDUAL.

Exactly as for the definition of the classical Multiple Correspondence Analysis, we introduce the Chi-2 distance simultaneously for the row and column profiles and we weight the individuals and the modalities proportionately to the frequency of each one. As we saw before, that is equivalent to do a correction of the complete disjunctive table, and to define:

$$d_{ij}^{c} = \frac{d_{ij}}{\sqrt{K\, d_{.j}}}.$$

where $K = \sum_{j=1}^{M} d_{ij}$ et $d_{.j} = \sum_{i=1}^{N} d_{ij}$.

$K$ is the number of questions and term $d_{.j}$ represents the number of persons who chose modality $j$.

When adjusted thusly, the table is called $D^c$ (Corrected Disjunctive Table).

We then consider a Kohonen network, and associate with each unit $u$ a code vector $C_u$ that is comprised of $(M + N)$ components, with the first $M$ components evolving in the space for individuals (represented by the rows of $D^c$) and the last $N$ components in the space for modalities (represented by the columns of $D^c$). We denote

$$C_u = (C_M, C_N)_u = (C_{M,u}, C_{N,u})$$

to put in evidence the structure of the code-vector $C_u$. The Kohonen algorithm lends itself to a double learning process. At each step, we alternatively draw a $D^c$ row (i.e. an individual $i$), or a $D^c$ column (i.e. a modality $j$).

When we draw an individual $i$, we associate the modality $j(i)$ defined by

$$j(i) = Arg \max_j d^c_{ij} = Arg \max_j \frac{d_{ij}}{\sqrt{Kd_{.j}}}$$

that maximizes the coefficient $d^c_{ij}$, i.e. the rarest modality out of all the corresponding ones in the total population. We then create an extended individual vector $X = (i, j(i)) = (X_M, X_N)$, of dimension $(M + N)$. See Fig. 2. Subsequently, we look for the closest of all the code vectors, in terms of the Euclidean distance restricted to the first $M$ components. Note $u_0$ the winning unit. Next we move the code vector of the unit $u_0$ and its neighbors closer to the extended vector $X = (i, j(i))$, as per the customary Kohonen law. Let us write down the formal definition :

$$\begin{cases} u_0 = Arg \min_u \left\| X_M - C_{M,u} \right\| \\ C_u^{new} = C_u^{old} + \varepsilon\,\sigma(u,u_0)\,(X - C_u^{old}) \end{cases}$$

where $\varepsilon$ is the adaptation parameter (positive, decreasing with time), and $\sigma$ is the neighborhood function, such that $\sigma(u, u_0)= 1$ if $u$ and $u_0$ are neighbour in the Kohonen network, and = 0 if not.

The reason to associate a row and a column in such a way is to keep the individual-modality associations which are realized in classical MCA by the fact that the principal axes of both diagonalizations are strongly related.

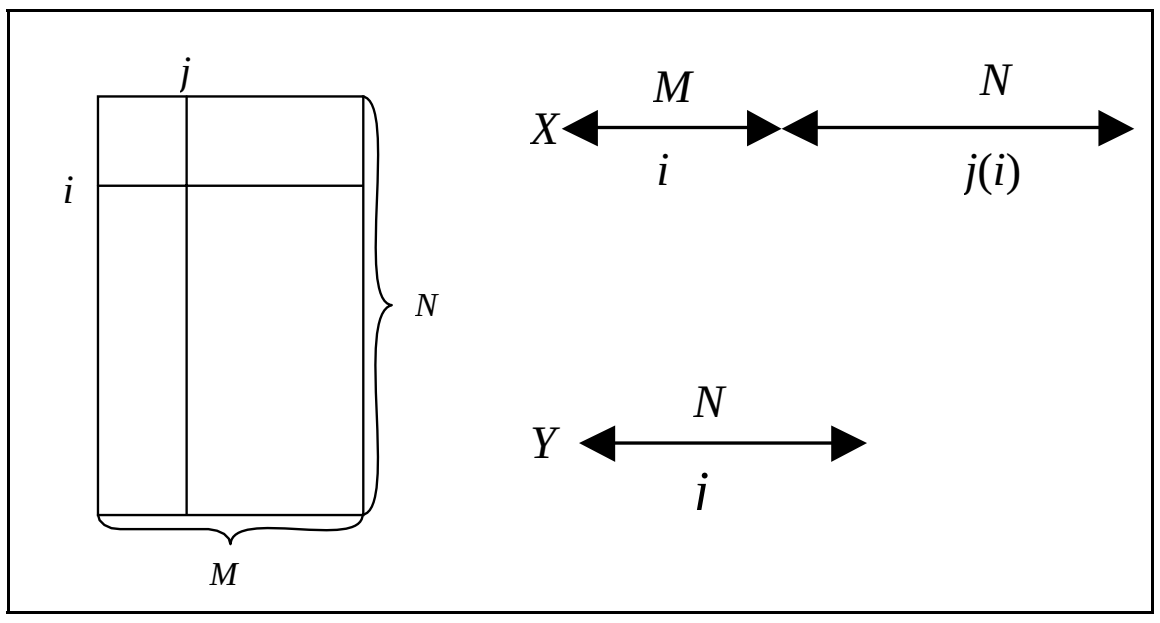

*Fig 2 : The matrix $D^c$, vectors X and Y*

When we draw a modality *j* with dimension *N* (a column of $D^c$), we do not associate an individual with it. Indeed, by construction, there are many equally placed individuals, and this would be an arbitrary choice. We then seek the code vector that is the closest, in terms of the Euclidean distance restricted to the last *N* components. Let $v_0$ be the winning unit. We then move the last *N* components of the winning code vector associated to $v_0$ and its neighbors closer to the corresponding components of the modality vector *j*, without modifying the first *M* components. For simplicity let us denote by *Y* (see Fig. 2) the *N*-column vector corresponding to modality *j*. This step can be written :

$$\begin{cases} v_0 = Arg\min_u \left\| Y - C_{N,u} \right\| \\ C_{N,u}^{new} = C_{N,u}^{old} + \varepsilon \quad \sigma(u, v_0) \quad (Y - C_{N,u}^{old}) \end{cases}$$

while the first *M* components are not modified.

This two-steps computation carries out a Kohonen classification of individuals (represented by the rows of the corrected Table $D^c$), plus a classification of modalities, maintaining all the while the associations of both individuals and modalities.

After convergence, the individuals and the modalities are classified into Kohonen classes. “Neighboring” individuals or modalities are classified in the same class or in neighboring classes. We call this algorithm just defined KDISJ. The name KDISJ is for Kohonen Disjunctive analysis. Its computing time is small, (the number of iterations is about 15 times the size of the database). In the following, we present two real data examples.

As we saw before, one of the main interest of the simultaneous classification of individuals and modalities is that it is possible to control the locations of the modalities by studying the distributions of the individuals that are classified in the same class and by computing the so-called deviations.

Let us recall what is a deviation for a modality *j*. The deviation for a modality *j* (shared by $d_{.j}$ individuals) and for a class *k* (with $n_k$, individuals) can be calculated as the difference between the number of individuals who possess this modality and belong to the class *k* and the "theoretical " number $d_{.j}\ n_k\ /\ N$ which would correspond to a distribution of the modality *j* in the class *k* that matches its distribution throughout the total population.

If the modalities are well located, all deviations have to be positive. In fact, when a deviation which is computed for modality *j* and class *k*, is positive, it means that class *k* gathers more individuals who share the modality *j* than it would do if the modalities were randomly distributed over the classes.

In the next section, we present a comparison between the KDISJ method and the classical MCA. To fulfil this comparison, we have to compare two classifications that are achieved in two different ways : one is the result of KDISJ, the other is achieved by classifying the individuals and the modalities after having computed their new coordinates resulting of the Multiple Correspondence Analysis.

## 3. Examples

### 3.1 Part-time employees

The data are extracted from a large INSEE 1998-1999 *Timetables* survey. See Letrémy et al, 2002, for the complete study (in French), Cottrell and Letrémy, 2003 for a first study about part-time employees. We only consider 207 part-time employees working on either an open-ended or a fixed term contract. They are described by 8 qualitative variables and 23 modalities according to the following table :

| **Heading** | **Modalities** | **Name** |
|---|---|---|
| **Type of employment contract** | Open-ended / fixed term contract | OEC,FTC |
| **Age** | <25, [25, 40], [40,50], ≥50 | AGE1,AGE2,AGE3,AGE4 |
| **Daily working schedules** | Identical, Posted, Variable | HOR1,HOR2,HOR3 |
| **Saturday shifts** | Never, sometimes, usually | SAT1,SAT2,SAT3 |
| **Sunday shifts** | Never, sometimes, usually | SUN1,SUN2,SUN3 |
| **Able to take time off** | Yes, yes under certain conditions, no | ABS1,ABS2,ABS3 |
| **Part-time status chosen** | Yes, no | CHO1,CHO2 |
| **Possibility of carrying over working hours** | Not applicable, yes, no | REC0,REC1,REC2 |

Simple cross analysis of the variables shows that the OEC contracts represent 83.57 % of all the population, while forced (and therefore involuntary) part-time work (CHO2) accounts for 46%. The goal is then to simultaneously represent all the modalities and the individuals, by realizing a classification of the (207 + 23) items. The KDISJ algorithm provides this classification, and as it is built by a SOM technique, the main associations and proximity are visible on the map.

On the map (Fig. 3), a 5 by 5 grid, using the KDISJ algorithm, we display findings from a simultaneous classification of individuals and variables. To simplify the representation, we have in each case displayed the modalities, the number of individuals (and not the complete list), the number of persons working on a chosen or not chosen part-time work. The starred units correspond to the classes where the number of fixed term contract (FTC) is greater than the mean value in the whole population (16.43 %).

We can see immediately that the bottom of the map correspond to the unpleasant working conditions (involuntary part-time status, fixed term contract, Saturday and Sunday shift, no chance to take any time off, etc... They are the youngest persons. The relatively favorable situations are displayed in the center of the map (with identical daily working schedules, no work on Saturday, nor Sunday, chance to take time off, chosen part-time status, open end contract, etc...).

| | | | | |
|---|---|---|---|---|
| SUN2<br>15 (13, 2) | HOR2<br>5 (4, 1) | REC1<br>16 (15, 1) | 8 (8, 0) | ABS2<br>13 (9, 4) |
| 8 (3, 5) | SAT2<br>3 (1, 2) | 2 (2, 0) | AGE3<br>13 (8, 5) | 4 (2, 2) |
| HOR3<br>REC2<br>12 (3, 9) | 6 (4, 2) | AGE4<br>17 (13, 4) | OEC<br>CHO1<br>3 (3, 0) | AGE2<br>16 (9, 7) |
| | AGE1<br>*<br>6 (2, 4) | 10 (1, 9) | HOR1,<br>SAT1,<br>SUN1,<br>ABS1,REC0 | *<br>5 (1, 4) |
| SUN3<br>13 (9, 4) | | SAT3<br>ABS3<br>11 (0, 11) | CHO2<br>*<br>3 (0, 3) | FTC<br>*<br>18 (1, 17) |

*Fig 3: The map of modalities and individuals.*

As pointed out before, it is possible to control the good position of the modalities with respect to the individuals, by computing the deviations. They are all positive, except for the modality ABS3, which should be classified in the class just above.

If we want to reduce the number of classes, it is possible to carry through a clustering of the 25 code vectors, by using any method,. For example, we use a one-dimensional Kohonen algorithm and group the 25 classes into 7 macro-classes, easier to describe, if the goal is to build a typology of all the individuals together with their modalities. The one-dimensional map classifies the contents along a decreasing scale from light gray (good conditions of work) to very dark gray (bad conditions of work).

The next figure (Fig. 4) summarizes this clustering into 7 macro-classes. We display the contents of each of them. For example, macro-class 1 gathers together classes 1, 2, 6 and 7 (very light gray in figure 3), and so on.

| | | | | | | |
|---|---|---|---|---|---|---|
| SUN2 | REC1 | AGE2 | HOR1 | SAT3 | HOR3 | SUN3 |
| HOR2 | ABS2 | AGE4 | SAT1 | ABS3 | REC2 | |
| SAT2 | AGE3 | OEC | SUN1 | CHO2 | AGE1 | |
| | | CHO1 | ABS1 | FTC | | |
| | | | REC0 | | | |
| 31 ind | 52 ind | 46 ind | 15 ind | 32 ind | 18 ind | 13 ind |

*Fig. 4: Seven organized macro-classes, after clustering the 25 classes.*

To compare these outputs with the results of a classical MCA, we represent in figure 5 the first projection (axes 1 and 2) of all individuals and modalities.

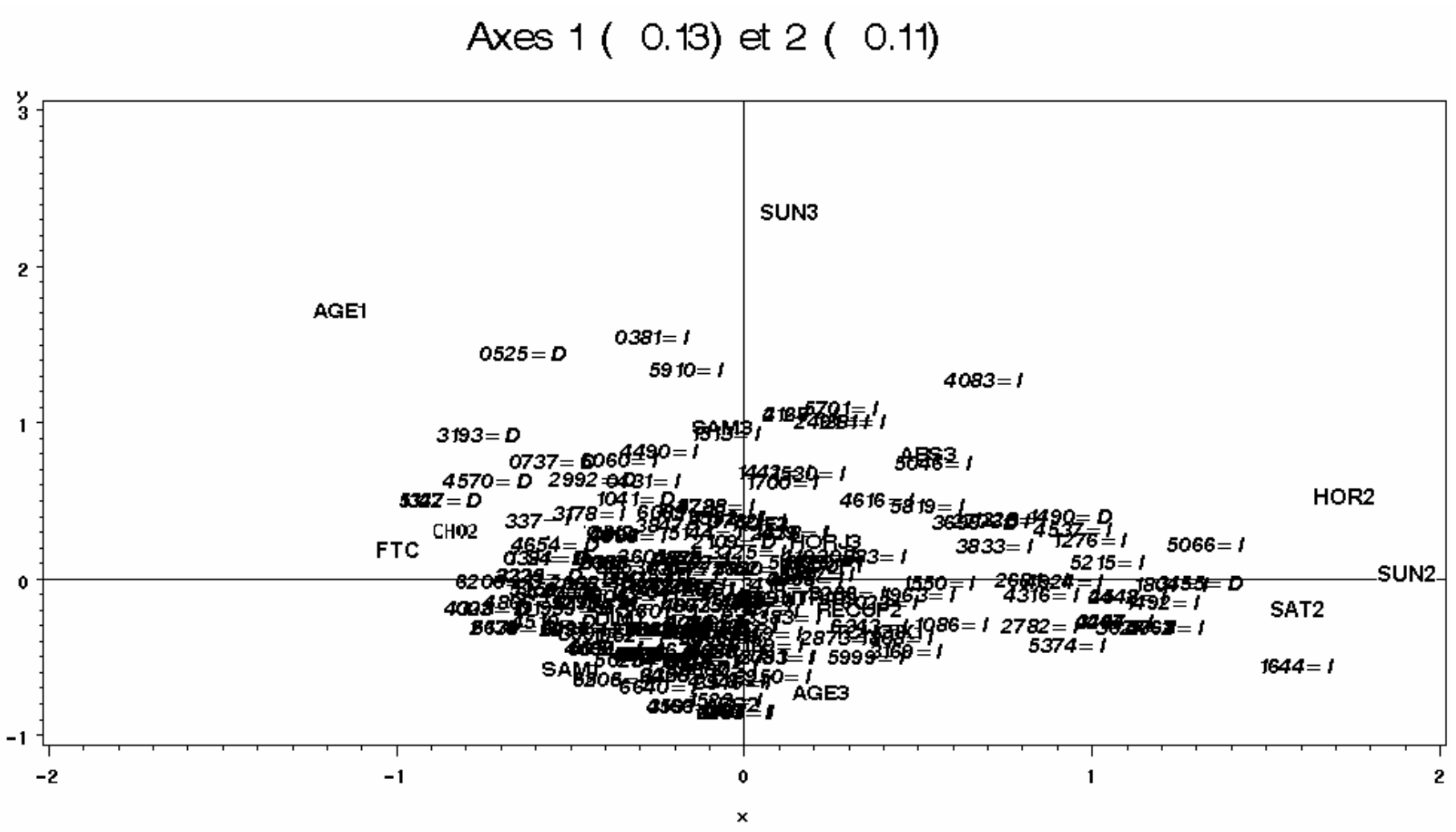

*Fig. 5 : MCA projection on axes 1 and 2 (24% of variance): individuals and modalities.*

We observe at once that the graphical representation is not readable, only some modalities can be seen, because they are on the borders of the figure. To improve the display, we represent in Fig. 6 the modalities alone. We can observe the same associations as in Fig. 3: AGE1 (the youngest people) is near bad conditions of work (FTC, CHO2, SAT3, SUN3, ABS3), and so on.

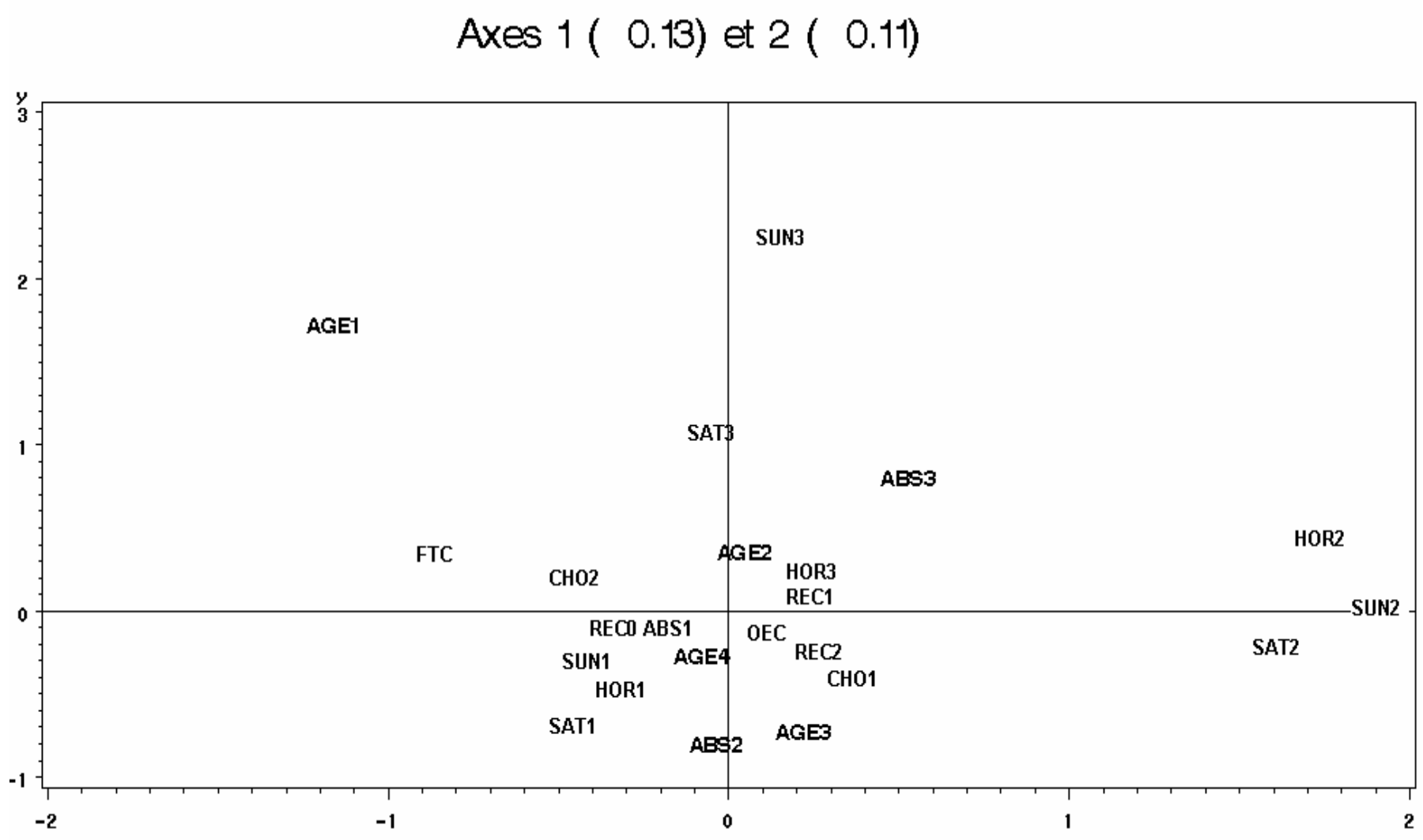

*Fig. 6 : MCA projection on axes 1 and 2 (24% of variance): modalities.*

One can try to simplify the representation by achieving a classification after the MCA transformation. Actually, after MCA, each individual and each modality is represented by new real-valued factorial coordinates on the $K$ - $M$ principal axes. In our example, there are 23 - 8 = 15 coordinates. So it is possible to simultaneously group the individuals and the modalities (for example with an Ascending Hierarchical Classification, AHC) into 25 classes (to get the same number of classes), and to compare the results with the classification that we got after using KDISJ.

Then one can compute the deviations according to the definition: they are also almost all positive except 2 (HOR3 and CHO2 have negative deviations). So from the point of view of classification, the result is not bad, even if it is not so good as the KDISJ classification, that has only one negative deviation. But the main drawback of this classification is that there is no neighboring relations between classes. So we lost the visualization properties of figure 3.

If a Kohonen classification is used on these transformed factorial coordinates, the results are similar, only two deviations are negative, and this once it is possible to display the classes on a map where the neighborhoods are apparent.

So the possible options are

| Method | Classification | Negative deviations | Visualization |
|---|---|---|---|
| KDISJ :<br>Correction of *D* + Kohonen algorithm | Yes | 1 | Good |
| MCA | Not | | Bad |
| MCA + AHC | Yes | 2 | Bad |
| MCA + Kohonen algorithm | Yes | 2 | Good |

Note that KDISJ save the computing time necessary to compute as well the eigenvalues and eigenvectors of matrix $D^{c}{}'D^{c}$ as the new coordinates of the individuals and modalities. That is particularly useful for large databases. It provides very quickly a good classification with nice visualization properties.

## 3.2 Recurring unemployed workers

This second example will be presented with less details for lack of space. The initial data is the complete register of the unemployed held by the ANPE; information on unemployment benefits and compensations in latest job is added from the data collected by UNEDIC. The studied period goes from July 1993 to August 1996. The population is constituted of all the unemployed who were looking for a job at the beginning of this period, or who became unemployed later (but before the end of August 1996); at the end of the period they are either unemployed or their status has changed in some way. In a previous study, (Cottrell and Gaubert, 2000), we use a 1% sample of the unemployed registered in the administrative region of Ile-de-France (Paris and suburbs) having two or more spells of unemployment (590 000 individuals on a population of more than 2 millions 167 000), those named recurring unemployed.

To illustrate the use of the KDISJ algorithm, we restrict ourselves in this paper to a sample of 204 individuals, described only by 8 qualitative attributes and 28 modalities. They are presented in the following table.

| | |
|---|---|
| **AGE** | Sub-categories of age: <25, 25-35, 35-45, 45-55, >55 |
| **BEN** | Daily benefits: <60 F, 60-100, 100-150, >150 |
| **EDU** | Level of education: > bac (post secondary school level),<br>bac level (secondary school completed),<br><bac (secondary school not completed) |
| **DUR** | Cumulated duration of unemployment: <12 months, 12-24, >24 |
| **OCC** | Monthly hours of occasional work OW : 0, 0-39, 39-78, 78-117, >117 |
| **POCC** | Proportion of cumulated duration of unemployment doing OW:<br>0, 0-0.1, 0.1-0.3, >0.3 |
| **EXIT** | Types of exits from unemployment (2 categories detailed below) |
| **REG** | Reasons for unemployment (2 categories detailed below) |

The different types of *exits from unemployment* have been grouped in 2 categories:

1. job found by the individual himself or with the help of ANPE services;
2. training program, withdrawal from the labor market, other exit.

Similarly, the causes of *registration at the ANPE* have been coded in 2 modalities:

1. lay off, end of fixed-term contracts, voluntary quit;
2. first job search.

Some of these categorical variables have an inherent order, in fact they were grouped into classes by the authors of the survey, to be able to compare with other studies. But in any case, the most important variables are the variables EXIT and REG, that are essentially categorical, and the goal is to deal with all variables simultaneously.

For this example, we use a one dimensional Kohonen map, with 6 units. The choice of a one-dimensional map is very useful in this study, since it allows to organize the situations in a very easily interpretable order, and to assign a mark (or score) to each class. We train it with the KDISJ algorithm, to simultaneously classify the modalities and the individuals. We get the following map, where the number at the bottom is the number of individuals in each class.

| | | | | | |
|---|---|---|---|---|---|
| OCC4<br>OCC5<br>POCC3,<br>POCC4 | OCC2<br>BEN2<br>OCC3<br>POCC2 | AGE5<br>AGE4<br>DUR3 | AGE3<br>BEN3<br>BEN4<br>EDU3 | AGE2, DUR2<br>OCC1, EXIT2<br>REG1 | AGE1, BEN1<br>EDU1, EDU2<br>DUR1,POCC1<br>EXIT1, REG2 |
| 45 | 21 | 22 | 40 | 34 | 42 |

*Fig 7: One-dimensional Kohonen map, modalities and individuals.*

We immediately see that the classes are well organized, according to the proportion of cumulated duration of occasional work, from left to right, according to the age. We see that contrarily to a very common idea, the duration of the unemployment is not associated to the daily benefit, but with the age.

Class 1 and 2 comprise people highly involved in occasional work. The benefits obtained are slightly above the average. Class 3 and 4 are characterized by the very long duration of unemployment. People are older than the average with benefits from unemployment slightly above the average. Most of them are not exerting any occasional work. Class 5 is constituted with young people, with no occasional work, and a short seniority in past employment. Their situation seems a good illustration of a typical trajectory with successive periods of fixed-term contracts and unemployment. Class 6 is made up of very young people, having a very short seniority, a duration of unemployment lower than the whole population mean, an average amount of unemployment benefits close to 0. Most of them are still looking for their first job. They leave unemployment by finding a job. We retrieve in a very easy way the main conclusions of the full study, (Cottrell and Gaubert, 2000).

## 4. Conclusion

It is possible to extend the basic SOM techniques to many other frames than the usual ones, like the study of individuals described by qualitative variables. Il is also very important to keep in mind that classical and Kohonen-based methodologies can be mixed. For example, it is possible to previously achieve a classification of the observations based on quantitative variables, and then to use the class number as a new qualitative variable which is added to the other ones, to apply the KDISJ algorithm. Conversely the KDISJ algorithm can be applied to the qualitative variables to get a partition of the data, to look for a specific model in each class. It is each day more and more evident that the Kohonen-based methods are a part of the numerous tools that the statistician have at his disposal to analyze, represent, visualize data.

**Acknowledgement**s: We thank ANPE and UNEDIC for the permission to use the data.